\documentclass[a4paper,twoside,11pt]{amsart}
\usepackage{fullpage}
\usepackage{amsmath} \usepackage{amssymb} \usepackage{amsopn}
\usepackage{amsthm} \usepackage{amsfonts} \usepackage{mathrsfs}

\usepackage{german}

\usepackage[pdfpagelabels,pdftex]{hyperref}
\usepackage[dvips]{graphicx} \usepackage[usenames]{color}
\usepackage[all]{xy}

\usepackage{setspace}
\usepackage{amssymb}
\usepackage{euscript}
\usepackage{cite}
\usepackage[all]{xy}
\usepackage{tabularx}

\theoremstyle{plain}
\usepackage{amsfonts}
\usepackage{graphicx}
\usepackage{amsmath}
\usepackage{amssymb}

\hypersetup{
  pdftitle={},
  pdfauthor={},
  pdfsubject={},
  pdfkeywords={},
  colorlinks=false,
  breaklinks=true,
  bookmarksopen=true,
  bookmarksnumbered=true,
  pdfpagemode=UseOutlines,
  plainpages=false}

\newcommand{\bG}{{\mathbb G}}

\newcommand{\bQ}{{\mathbb Q}}
\newcommand{\bR}{{\mathbb R}}

\newcommand{\bZ}{{\mathbb Z}}

\newcommand{\cB}{{\mathscr B}}

\newcommand{\fH}{{\mathfrak H}}

\newcommand{\fo}{{\mathfrak o}}

\DeclareMathOperator{\GL}{GL}

\newtheorem{thm}{Theorem}[section]

\newtheorem{prop}[thm]{Proposition}

\newtheorem{cor}[thm]{Corollary}

\newtheorem{lm}[thm]{Lemma}

\theoremstyle{definition}

\theoremstyle{remark}
\newtheorem{rem}[thm]{Remark}

\newenvironment{pro*}[1][Proof]{{\it{#1:}} }{}

\def\abs#1{\left|#1\right|}

\newcommand\minn{\mathop{ \rm min}}

\newcommand\bk{\mathop{ \bar{k} }}

\newcommand\rar{ \rightarrow }
\newcommand\tar{ \twoheadrightarrow }

\DeclareMathOperator{\disjun}{\dot{\cup}}
\DeclareMathOperator{\card}{{\rm card}}
\newcommand{\sm}{{\,\smallsetminus\,}}

\newcounter{absatzcounter}[section]
\numberwithin{equation}{section}

\begin{document}

\title{On a finite type property of certain affine Deligne-Lusztig varieties}
\author{Alexander Ivanov}
\maketitle

\section{Introduction}

Let $k$ be a field with $q$ elements, and let $\bar{k}$ be an algebraic closure of $k$. Let $\sigma$ denote the Frobenius morphism of $\bar{k}/k$. Put $F = k((t))$ and $L = \bar{k}((t))$. We extend $\sigma$ to the Frobenius morphism of $L/F$ by setting $\sigma(t) = t$. Write $\mathfrak{o} = \bar{k}[[t]]$ for the valuation ring of $L$.

Let $G$ be a split connected reductive group over $k$ and let $A$ be a split maximal torus in $G$. For $b \in G(L)$, let $\nu_b \in X_{\ast}(A)_{\bQ}$ denote the Newton point of $b$.
Write $W$ and $\tilde{W} = X_{\ast}(A) \rtimes W$ for the finite and the extended affine Weyl groups attached to $A$. Fix a Borel subgroup $B$ containing $A$ and let $I$ be the preimage of $B(\bk)$ under the projection $G(\fo) \twoheadrightarrow G(\bk)$. Then $I$ is an Iwahori subgroup in $G(L)$. Let $X = G(L)/I$ be the affine flag manifold. The group $G(L)$ acts on $X$ by left translation. The Bruhat decomposition implies that $G(L)$ is the union of the double cosets $IwI$, where $w \in \tilde{W}$. Following \cite{Ra}, the affine Deligne-Lusztig variety $X_w(b)$ attached to $b \in G(L)$ and $w \in \tilde{W}$ is the locally closed subset of $X$, endowed with its reduced induced sub-Ind-scheme structure, defined by
\begin{equation*}
 X_w(b) = \{ xI \in G(L)/I \colon x^{-1}b \sigma(x) \in IwI\}.
\end{equation*}

\noindent Then $X_w(b)$ is locally of finite type, which follows from \cite{RZ} theorem 1.4, but in general
not of finite type. But if $G$ is of adjoint type and $b$ is
superbasic, i.e. $\nu_b$ is central and $b$ lies in no proper Levi
subgroup of $G$, then $X_w(b)$ is of finite type. This was proven by
Viehmann in \cite{Vi}, among other things. In this short note, we
give an alternative proof. Since superbasic elements occur only for
type $A_n$, we are reduced to show the following (compare also
\cite{Go:Uebersicht}, 4.13).

\begin{prop}[\cite{Vi}]\label{prop:ftsimadjXwb} Let $G = PGL_n$, $b \in G(L)$ superbasic and $w \in \tilde{W}$. Then $X_w(b)$ is of finite type. \end{prop}

\noindent To prove this, we will show that there are only finitely
many $v \in \tilde{W}$, such that $IvI/I \cap X_w(b) \neq
\emptyset$. The main ingredient in the proof is proposition
\ref{prop:LinFunctBeschr}. Let $\beta$ denote the automorphism of
$G(L)$ defined by conjugation with $b$. It induces an automorphism
of the affine Weyl group $W_a \subseteq \tilde{W}$. Then proposition
\ref{prop:LinFunctBeschr} gives a lower bound for
$\ell(\beta(v)v^{-1})$ by an expression, which is linear in
$\ell(v)$. Such an estimate is a special case of a result of
Rousseau (\cite{Ro} proposition 2.3; compare also Rapoport-Zink
\cite{RZ} theorem 1.4), who considers a building $\cB$ (and in
particular the Bruhat-Tits building of a reductive algebraic group),
and shows that if $\beta$ is an automorphism of $\cB$, having a
non-empty fixed point set $\cB^{\langle \beta \rangle }$, and $d$
denotes the distance function on $\cB$, then there is a constant $c
> 0$, depending only on the geometry of $\cB$ (not on $\beta$!), such that for all $x \in \cB$, $\frac{d(x, \beta(x))}{d(x,
\cB^{\langle \beta \rangle })} > c$. Our proof differs from \cite{Ro}: the estimate is
formulated in terms of the affine Weyl group; the proof is explicit, works only in the special sitiation, and gives some
quantitative information on the length of involved Weyl group
elements. This gives an estimate of the dimension of $X_w(b)$, which
is, however, rather weak. For $n = 2$ it gives $\dim X_w(b) \leq
\frac{1}{2} \ell(w) + 2$ for all $w \in \tilde{W}$, which is sharp
up to a constant. But already for $n=3$, an improved version of the
estimate from proposition \ref{prop:LinFunctBeschr} gives $\dim
X_w(b) \leq \frac{2}{3}\ell(w) + 7$ for all $w \in \tilde{W}$, which
is weaker than results proven in \cite{GH}. The proof of
\ref{prop:ftsimadjXwb} is given in section
\ref{sec:FiniteTypeProperty}.


In sections \ref{sec:A propertyOfTheLength} and \ref{sec:HeckeAlgebra}, we show a property of the Hecke algebra of a Coxeter group, which we need in our proof. Let $(W,S)$ be a Coxeter system. Let $\fH = \fH(W)$ be the Hecke algebra associated with (W,S), i.e. a $\bZ[v,v^{-1}]$-algebra generated by elements $T_s$ with $s \in S$, with certain relations (see section \ref{sec:HeckeAlgebra} or \cite{Lu} 3.2). If $x = s_1 s_2 \dots s_n \in W$ is a reduced expression, we write $T_x = T_{s_1} \dots T_{s_n}$. The set $\{T_x \colon x \in W\}$ is a basis of $\fH(W)$ as a free $\bZ[v,v^{-1}]$-module, and hence for any two elements $x,y \in W$, we can write
\[T_x T_y = \sum_{w \in W} r_w^{x,y} T_w, \]
\noindent with $r_w^{x,y} \in \bZ[v,v^{-1}]$. One can ask, how the set
\[D(x,y) = \{ w \in W \colon r_w^{x,y} \neq 0 \} \]
\noindent looks like (this is one of the questions studied in \cite{Ha}). We give a proof of the following property of it.

\begin{prop}\label{prop:HeckeAlgMult} Let $x,y \in W$. Then $r_w^{x,y} = 0$ unless $\ell(xy) \leq \ell(w) \leq \ell(x) + \ell(y)$. \end{prop}

\noindent The second inequality is trivial, the first follows almost immediately from the proposition \ref{prop:miracleprop} below. In particular, if $G$ is  a split connected reductive group over $k$ and $\tilde{W}$ the extended affine Weyl group of $G$, this can be interpreted as follows: the product of two Schubert cells in the affine flag manifold of $G$ attached to $x,y \in \tilde{W}$ is the union of Schubert cells of dimensions $\geq \ell(xy)$ (the same is also true in the situation of a finite root system).

To prove proposition \ref{prop:HeckeAlgMult} we need (a weaker version of) a result on general Coxeter groups, proven for example in (\cite{BB} lemma 2.2.10) or (\cite{Ha} lemma 5.6). We reprove it in section \ref{sec:A propertyOfTheLength}, omitting the direct use of the strong exchange property, in contrast to both references.

\begin{prop}\label{prop:miracleprop}
Let $(W,S)$ be a Coxeter system. Let $s \in S$ be a simple reflection and $x,y \in W$, such that $\ell(xs) > \ell(x)$ and $\ell(sy) > \ell(y)$. Then $\ell(xsy) > \ell(xy)$.
\end{prop}

\textbf{Acknowledgments.} I am very grateful to Ulrich G\"{o}rtz for helpful discussions on this subject and for pointing out to me some of the references in this note. Also, I am grateful to Juan Cervi\~{n}o for some remarks and the interest in my work.

\section{A property of the length in a Coxeter group}\label{sec:A propertyOfTheLength}
In this seciton $(W,S)$ denotes a Coxeter system. We will prove proposition \ref{prop:miracleprop}. Here is an immediate corollary from it:


\begin{cor} Under the assumptions of the proposition \ref{prop:miracleprop}, assume additionally $\ell(xy) = \ell(x) + \ell(y)$. Then $\ell(xsy) = \ell(x) + \ell(y) + 1$.
\end{cor}


%

One can attach to $(W,S)$ a geometric realization and a root system. Details can be found in (\cite{Bo}, Chap V, \S4), \cite{St}, or \cite{De}. We briefly recall the construction in the section \ref{sec:roots} following the last reference.

\subsection{Root system attached to $(W,S)$} \label{sec:roots}
For $s,t \in S$, denote by $m_{st}$ the order of $st$ in $W$. Let $E$ be the real vector space with the basis $\{e_s : s \in S\}$. Then $W$ determines the symmetric bilinear form $(,)$ on $E$ defined by:
\[(e_s,e_{s^{\prime}}) = -\cos(\pi/m_{ss^{\prime}}), \text{ for } s,s^{\prime} \in S\]

\noindent (if $m_{ss^{\prime}} = \infty$, then $(e_s, e_{s^{\prime}}) := - 1$). We have $(e_s,e_s) = 1$. There exists a unique representation
\[\sigma \colon W \rar \GL(E),\]
subjected to the condition that $\sigma(s)(e) = e  - 2(e,e_s)e_s$ for all $s \in S$ and all $e \in E$. This representation is faithful and we call it the \emph{geometric realization} of $(W,S)$. We left $\sigma$ out of the notation and write $x.e$ instead of $\sigma(x)(e)$. Further, $(,)$ is $W$-invariant. Now let
\[\Phi := \{ x.e_s \colon x \in W, s \in S \}.\]
be the set of roots. We have $(\alpha,\alpha) = 1$ for all $\alpha \in \Phi$.

\begin{prop} (\cite{Bo}, Chap V, \S4, ex. 8)
Let $\alpha \in \Phi$. Then $\alpha = \sum_{s \in S} a_s e_s$, where either all $a_s \geq 0$ or all $a_s \leq 0$.
\end{prop}

\noindent This can be proved by induction on the length. It allows us to define the (disjoint) partition $\Phi =: \Phi^+ \disjun \Phi^-$, where
\[ \Phi^+ := \{\alpha \in \Phi \colon \alpha = \sum_{s \in S} a_s e_s \text{ with } a_s \geq 0 \text{ for all } s \in S \} \]
and $\Phi^- := - \Phi^+$. For a root $\alpha$, we write $\alpha > 0$ if $\alpha \in \Phi^+$ and $\alpha < 0$ if $\alpha \in \Phi^-$. For $x \in W$ set:
\begin{eqnarray*} \Phi_x &:=& \{\alpha \in \Phi \colon \alpha > 0, x.\alpha < 0 \} \quad \text{ and } \\
 \Phi^-_x &:=& - \Phi_x.
\end{eqnarray*}

We have the following fundamental facts:
\begin{prop}\label{prop:PropOfRoots} (\cite{St}, \S 1) Let $s \in S$, $x \in W$.
\begin{itemize}
\item[(i)] $\Phi_s = \{e_s\}$, or equivalently $s.(\Phi^+ \sm \{e_s\}) = (\Phi^+ \sm \{e_s\})$.
\item[(ii)] $\ell(xs) > \ell(x) \Leftrightarrow e_s \not\in \Phi_x$.
\item[(iii)] $\ell(x) = \card(\Phi_x)$.
\end{itemize}
\end{prop}

\noindent There is a partial order on $E$: $\alpha \geq \beta$ if and only if $\alpha - \beta$ is a non-negative linear combination of positive roots.

\subsection{Further properties of roots}

In this section we prove two lemmas needed later.

\begin{lm}\label{lm:formulas} Let $x, y \in W$. Then:
\begin{itemize}
\item[(i)] $\Phi_{x^{-1}} = x.\Phi^-_x$
\item[(ii)] $\Phi_{xy} = (\Phi_y \sm y^{-1}.\Phi^-_x) \disjun (y^{-1}.\Phi_x \sm \Phi^-_y)$.
\end{itemize}
\end{lm}
\begin{proof} (i): a substitution $\beta = x\alpha$ gives:
\[ x.\Phi^-_x = x.\{\alpha \colon \alpha < 0, x.\alpha > 0\} = \{\beta \colon \beta > 0, x^{-1}.\beta < 0 \} = \Phi_{x^{-1}}.\]
(ii): Write $\Phi_{xy}$ as the disjoint union
\begin{eqnarray*} \Phi_{xy} = \{\alpha \colon \alpha > 0, xy.\alpha < 0\} &=& \{\alpha \colon \alpha > 0, y.\alpha < 0, xy.\alpha < 0\} \disjun \{\alpha \colon \alpha > 0, y.\alpha > 0, xy.\alpha < 0\} \\ &=& \{ \alpha \in \Phi_y \colon xy.\alpha < 0 \} \disjun y^{-1}.\{\beta \colon \beta > 0, x.\beta < 0, y^{-1}.\beta > 0 \}.
\end{eqnarray*}

Now,
\begin{eqnarray*} \{ \alpha \in \Phi_y \colon xy.\alpha < 0 \} &=& \Phi_y \sm y^{-1}.\Phi^-_x \quad \text{ and } \\
 y^{-1}.\{\beta \colon \beta > 0, x.\beta < 0, y^{-1}.\beta > 0 \} &=& y^{-1}.(\Phi_x \sm \Phi_{y^{-1}}) = y^{-1}.\Phi_x \sm \Phi^-_y,
\end{eqnarray*}
where the last equality is a consequence of (i).
\end{proof}

\noindent For $\beta \in \Phi^+$ and $s \in S$, we have $s.\beta = \beta - 2 (\beta, e_s) e_s$. I.e. either $s.\beta = \beta$, or $s.\beta > \beta$ or $s.\beta < \beta$.

\begin{lm}\label{lm:pairs} Let $s \in S$.
Let $x \in W$ with $\ell(xs) > \ell(x)$ and $\beta \in \Phi^+ \sm \{e_s\}$ with $s \beta \geq \beta$.
Then \[ s.\beta \in \Phi_x \Rightarrow \beta \in \Phi_x.\]
\end{lm}
\begin{proof} By assumption, we have $s.\beta - \beta = n e_s$ with $n \in \bR_{\geq 0}$. Assume $s \beta \in \Phi_x$, that is $xs.\beta = \sum_{t \in S} a_t e_t$ with $a_t \leq 0$. We have to show $\beta \in \Phi_x$, that is $x.\beta < 0$. But $x.\beta = xs.\beta - xs.\beta + x.\beta = xs.\beta - x.(s.\beta - \beta) = xs.\beta - x.(n e_s) = \sum_t a_t e_t - n(x.e_s)$. But by assumption and proposition \ref{prop:PropOfRoots}(ii) we have $x.e_s > 0$, i.e. $ - n(x.e_s) = \sum_t b_t e_t$ with $b_t \leq 0$. Therefore $x.\beta = \sum_t (a_t + b_t)e_t \in \Phi^-$.
\end{proof}


\subsection{Proof of proposition \ref{prop:miracleprop}}

Assume $x,y \in W$, $s \in S$ are given such that $\ell(xs) > \ell(x), \ell(sy) > \ell(y)$. This is equivalent to $e_s \not\in \Phi_x$, $e_s \not\in \Phi_{y^{-1}}$. Using lemma \ref{lm:formulas}(ii) we get:
\begin{eqnarray*}
\Phi_{xs} &=& s\Phi_x \disjun \{e_s\}, \\
\Phi_{sy} &=& \Phi_y \disjun \{y^{-1}.e_s \}.
\end{eqnarray*}
A further application of lemma \ref{lm:formulas}(ii) gives:
\begin{eqnarray*}
\Phi_{xy} &=& (\Phi_y \sm y^{-1}.\Phi^-_x) \disjun (y^{-1}.\Phi_x \sm \Phi^-_y), \\
\Phi_{xsy} &=& ((\Phi_y \disjun \{y^{-1}.e_s\}) \sm y^{-1}s.\Phi^-_x) \disjun (y^{-1}s.\Phi_x \sm (\Phi^-_y \disjun \{-y^{-1}.e_s\} )) \\
&=& (\Phi_y \sm y^{-1}s.\Phi^-_x) \disjun (y^{-1}s.\Phi_x \sm \Phi^-_y) \disjun \{ y^{-1}.e_s \},
\end{eqnarray*}
where the last equality follows from $e_s \not\in \Phi_x$. We set:
\begin{eqnarray*}
A_1 &:= \Phi_y \sm y^{-1}.\Phi^-_x \qquad A_2 &:= y^{-1}.\Phi_x \sm \Phi^-_y, \\
B_1 &:= \Phi_y \sm y^{-1}s.\Phi^-_x \qquad B_2 &:= y^{-1}s.\Phi_x \sm \Phi^-_y,
\end{eqnarray*}
i.e. $\Phi_{xy} = A_1 \disjun A_2$, $\Phi_{xsy} = B_1 \disjun B_2 \disjun \{ y^{-1}.e_s \}$. Hence we can write
\begin{eqnarray*}
\Phi_{xy} &=& (A_1 \cap B_1) \disjun (A_1 \sm B_1) \disjun (A_2 \cap B_2) \disjun (A_2 \sm B_2) \\
\Phi_{xsy} &=& (A_1 \cap B_1) \disjun (B_1 \sm A_1) \disjun (A_2 \cap B_2) \disjun (B_2 \sm A_2) \disjun \{ y^{-1}.e_s \}.
\end{eqnarray*}

We claim that $\card(A_1 \sm B_1) \leq \card(B_1 \sm A_1)$ and $\card(A_2 \sm B_2) \leq \card(B_2 \sm A_2)$. Since $\ell(w) = \card(\Phi_w)$ for any $w \in W$, the claim implies the assertion of the proposition.

Let us first proof that $\card(A_1 \sm B_1) \leq \card(B_1 \sm A_1)$. More precise, we claim that $\beta \mapsto y^{-1}sy.\beta$ defines an injection from $A_1 \sm B_1$ into $B_1 \sm A_1$. It is enough to show that if $\beta \in A_1 \sm B_1$, then $y^{-1}sy.\beta \in B_1 \sm A_1$. Thus let $\beta \in A_1 \sm B_1$ and set $\gamma := -y.\beta$. We have
\[ \beta \in A_1 \sm B_1 \Leftrightarrow \beta \in (\Phi_y \cap y^{-1}s.\Phi^-_x) \sm y^{-1}.\Phi^-_x \Leftrightarrow \gamma \in (y.\Phi^-_y \cap s.\Phi_x ) \sm \Phi_x = (\Phi_{y^{-1}} \cap s.\Phi_x ) \sm \Phi_x.\]
In particular $\gamma > 0$, $\gamma \neq e_s$ (since $\gamma \in \Phi_{y^{-1}}$) and $s.\gamma < \gamma$: otherwise we would have $s.\gamma \geq \gamma$ and $s\gamma \in \Phi_x$ would imply $\gamma \in \Phi_x$ by lemma \ref{lm:pairs}. But using lemma \ref{lm:pairs} again, we see that $s.\gamma \in \Phi_{y^{-1}}$, since $\gamma \in \Phi_{y^{-1}}$ and $s.\gamma < \gamma$. Therefore we obtain ($s.\gamma \in \Phi_{y^{-1}}, s.\gamma \in \Phi_x, s.(s.\gamma) \not\in \Phi_x$), i.e. $s.\gamma \in (\Phi_{y^{-1}} \cap \Phi_x) \sm s.\Phi_x$ and hence
\[y^{-1}sy.\beta = -y^{-1}s.\gamma \in y^{-1}.((\Phi^-_{y^{-1}} \cap \Phi^-_x) \sm s.\Phi^-_x) = (\Phi_y \cap y^{-1}.\Phi^-_x) \sm y^{-1}s.\Phi_x^- = B_1 \sm A_1.\]

Secondly, we have to prove that $\card(A_2 \sm B_2) \leq \card(B_2 \sm A_2)$. Analogously, we claim that $\beta \mapsto y^{-1}sy.\beta$ defines an injection from $A_2 \sm B_2$ into $B_2 \sm A_2$. Assume $\beta \in A_2 \sm B_2$. We have to prove $y^{-1}sy.\beta \in B_2 \sm A_2$. Set $\gamma := y.\beta$. Then
\[ \beta \in A_2 \sm B_2 \Leftrightarrow \beta \in y^{-1}.\Phi_x \sm (\Phi^-_y \cup y^{-1}s.\Phi_x) \Leftrightarrow \gamma \in \Phi_x \sm (\Phi_{y^{-1}} \cup s.\Phi_x).\]
In particular, $\gamma > 0$, $\gamma \neq e_s$ (since $\gamma \in \Phi_x$). As $\gamma \in \Phi_x$, $s.\gamma \not\in \Phi_x$, we obtain from lemma \ref{lm:pairs} $s.\gamma > \gamma$. This and lemma \ref{lm:pairs} assert $s.\gamma \not\in \Phi_{y^{-1}}$. Thus we obtain ($s.(s.\gamma) \in \Phi_x, s.\gamma \not\in \Phi_x, s.\gamma \not\in \Phi_{y^{-1}}$), i.e. $s.\gamma \in s.\Phi_x \sm (\Phi_x \cup \Phi_{y^{-1}})$. Hence
\[ y^{-1}sy.\beta = y^{-1}s.\gamma \in y^{-1}.(s.\Phi_x \sm (\Phi_x \cup \Phi_{y^{-1}})) = y^{-1}s.\Phi_x \sm (y^{-1}.\Phi_x \cup \Phi^-_y) = B_2 \sm A_2.\]
This finishes the proof. \qed

\section{Hecke Algebra}\label{sec:HeckeAlgebra}

In this section we prove proposition \ref{prop:HeckeAlgMult}.

\subsection{Definition} The Hecke algebra $\fH$ of $W$ with respect to any weight function $L \colon W \rar \bZ$, as defined in \cite{Lu} 3.2, is a $\bZ[v,v^{-1}]$-algebra generated by elements $T_s$ with $s \in S$, with relations
\begin{itemize}
 \item[(i)] $(T_s - v_s)(T_s + v_s^{-1}) = 0$,
\item[(ii)] $T_{s_1} T_{s_2} T_{s_1} \dots$ $= T_{s_2} T_{s_1} T_{s_2} \dots$,
\end{itemize}
for any $s,s_1,s_2 \in S$, where the number of factors on each side in the second line is equal to the order of $s_1 s_2$, and where $v_s := v^{L(s)}$. For an element $x \in W$ with reduced decomposition $x = s_1 s_2 \dots s_n$, we write as usual $T_x := T_{s_1} T_{s_2} \dots T_{s_n}$. Then $\fH$ is free $\bZ[v,v^{-1}]$-module with basis $\{T_w \colon w \in W\}$. As in the introduction, for any two elements $x,y \in W$, we can write
\[T_x T_y = \sum_{w \in W} r_w^{x,y} T_w, \]
\noindent with $r_w^{x,y} \in \bZ[v,v^{-1}]$.
Then proposition \ref{prop:HeckeAlgMult} states that if $r_w^{x,y} \neq 0$, then  $\ell(w) \geq \ell(xy)$. We prove it below.

\subsection{Proof of proposition \ref{prop:HeckeAlgMult}} \label{sec:ProofOfCorHeckeAlgMult}

Let $D^{\prime}(x,y)$ be defined by induction on $\ell(y)$ through
\[ D^{\prime}(x,y) := \begin{cases} \{x\} & \text{ if } y = 1, \\ D^{\prime}(xs,sy) & \text{ if } \ell(sy) < \ell(y) \text{ and } \ell(xs) > \ell(x), \\ D^{\prime}(xs, sy) \cup D^{\prime}(x,sy) & \text{ if } \ell(sy) < \ell(y) \text{ and } \ell(xs) < \ell(x). \end{cases} \]

Let $D(x,y) := \{w \in W \colon r_w^{x,y} \neq 0 \}$. Then, using the relations defining $\fH$, we obtain by induction on $\ell(y)$, that $D(x,y) \subseteq D^{\prime}(x,y)$. Let $m(x,y) := \min\{ \ell(w) \colon w \in D^{\prime}(x,y) \}$. We prove by induction on $\ell(y)$, that $m(x,y) \geq \ell(xy)$ (the converse inequality is easy to see, since $xy \in D^{\prime}(x,y)$). For $\ell(y) = 0$ there is nothing to do. Assume $\ell(y) > 0$ and let $s \in S$ be such that $\ell(y) > \ell(sy)$. If $\ell(x) < \ell(xs)$, then $D^{\prime}(x,y) = D^{\prime}(xs,sy)$ and we can use the induction hypothesis to see that $m(x,y) = m(xs,sy) \geq \ell(xssy) = \ell(xy)$. If $\ell(x) > \ell(xs)$, then $D^{\prime}(x,y) = D^{\prime}(xs,sy) \cup D^{\prime}(x,sy)$. By induction hypothesis we obtain: $m(xs,sy) \geq \ell(xssy) = \ell(xy)$, $m(x,sy) \geq \ell(xsy) > \ell(xy)$, where the last inequality follows from proposition \ref{prop:miracleprop}, applied to $xs,s,sy$. This finishes the proof.

\begin{rem} The corollary has the following geometric interpretation in the case of an affine root system. If $G$ is a split connected reductive group over a field $k$ and $I$ is an Iwahori subgroup of $G(k((\epsilon)))$, then for two elements $x,y \in \tilde{W}$ in the extended affine Weyl group of $G$, the product $IxIyI/I$ of the Schubert cells attached to $x,y$ is the exactly the union
\[IxIyI/I = \bigcup_{v \in D(x,y) } IvI/I. \]
Since the dimension of the Schubert cell associated to $v$ is equal to $\ell(v)$, the corollary shows that the dimensions of Schubert cells occurring in the decomposition of $IxIyI/I$ are $\geq \ell(xy)$. The same holds for a finite root system.
\end{rem}

\section{Finite type property of $X_w(b)$ with $b$ superbasic}\label{sec:FiniteTypeProperty}

In this section we prove proposition \ref{prop:ftsimadjXwb} from the introduction.

\subsection{Some preliminaries}\label{sec:someprelsXwb}
We use the notation from the introduction and we set $G = \GL_n$ with $n \geq 2$. We can assume that $A$ is the diagonal torus and $B$ the Borel subgroup of upper triangular matrices. Consider the $L$-vector space $L^n$ with standard basis $e_0, \dots, e_{n-1}$, on which $G(L)$ acts on the left. As in \cite{Vi}, we define $e_i$ for $i \in \bZ$ by $e_{i+n} := \epsilon e_i$. Since $b$ is superbasic, its Newton point will be of the form $\nu_b = (\frac{m}{n}, \dots, \frac{m}{n})$ with $m$ coprime to $n$. Without loss of generality, we choose $b$ to be the representative of its $\sigma$-conjugacy class, given by $e_i \mapsto e_{i+m}$.
The connected components of $X = G(L)/I$ are naturally indexed by $\pi_1(G) \cong \bZ$, the morphism mapping $gI$ to its connected component given by $gI \mapsto v_L(\det(g))$. Let $X_w(b)^i$ denote the intersection of the $i$-th connected component $X$ with $X_w(b)$. Let $b_1 \in \GL_n(L)$ be the element defined by $e_i \mapsto e_{i+1}$. Then $b = b_1^m$ and in particular $b_1$ commutes with $b$, and thus maps $X_w(b)^i$ isomorphically onto $X_w(b)^{i+1}$. Furthermore, $b$ commutes with $I$.

We denote by $s_1, \dots, s_{n-1}$ the reflections attached to the finite simple roots determined by the choice of $B$, numbered in an obvious way. They generate the finite Weyl group $W$ of $G$. Further let $s_0 := \epsilon^{\theta^{\vee}} s_{\theta} \in W_a \subset \tilde{W}$, where $W_a$ is the affine Weyl group, $\theta$ the longest finite positive root and $s_{\theta} \in W$ the associated reflection. Then $\{s_0, s_1, \dots, s_{n-1} \}$ is the distinguished set of Coxeter generators of $W_a$, given by the choice of $B$. Further, we extend the notation by setting $s_i = s_{i \mod n}$ for all $i \in \bZ$. Similarly, for any $\lambda = (\lambda_1, \dots, \lambda_n) \in X_{\ast}(A)$, and any $i \in \bZ$, we write $\lambda_i := \lambda_j$ where $0 < j \leq n$ is unique integer such that $i \equiv j \mod n$. Let $\beta$ denote the inner automorphism of $G(L)$ defined by $b$. It commutes with $\sigma$ and induces an (outer) isomorphism of $W_a$. An easy computation shows that  $\beta(s_i) = s_{i + m}$ and $\beta(\epsilon^{(\lambda_1, \dots, \lambda_n)}) = \epsilon^{(\lambda_{1 + m}, \dots, \lambda_{n + m})}$.

Since $X_w(b)$ for $PGL_n$ is isomorphic to a union of $n$ copies of $X_w(b)^0$ for $GL_n$ (an isomorphism being given by projection $GL_n \tar PGL_n$ and $b_1^i \colon X_w(b)^0 \stackrel{\sim}{\rar} X_w(b)^i$ as above), to prove proposition \ref{prop:ftsimadjXwb} it is enough to prove that $X_w(b)^0$ for $GL_n$ are of finite type. Furthermore, we can assume that $X_w(b)$ is non-empty, and in particular that $w$ and $b$ lie in the same connected component of $G/I$. Thus we can write $w = w_a b$ with $w_a \in W_a$.

\subsection{Proof of proposition \ref{prop:ftsimadjXwb}}

By the remark at the end of section \ref{sec:someprelsXwb}, we have to prove that $X_w(b)^0$ is of finite type. It is enough to show that there are only finitely many $v \in \tilde{W}$, with $IvI/I \cap X_w(b)^0 \neq \emptyset$. First of all, remark that such $v$ must necessarily satisfy $v \in W_a \subseteq \tilde{W}$, since the valuation of  its determinant must be $0$. If $v \in W_a$ is such that $gI \in IvI/I \cap X_w(b)^{0} \neq \emptyset$, then we have
\[ I w_a I b = IwI = I g^{-1} b \sigma(g) I = I g^{-1} \beta(\sigma(g)) I b \subseteq I g^{-1} I \sigma(\beta(v)) I b = Iv^{-1}I \beta(v) I b. \]

\noindent Thus $I w_a I \subseteq I v^{-1} I \beta(v) I$, i.e. $w_a \in D(v^{-1}, \beta(v))$, using the notation of section \ref{sec:ProofOfCorHeckeAlgMult}. Thus $\ell(w_a) \geq \ell(v^{-1} \beta(v)) = \ell(\beta(v^{-1})v)$ by proposition \ref{prop:HeckeAlgMult}. By proposition \ref{prop:LinFunctBeschr} applied to $v^{-1}$, there are only finitely many $v \in W_a$ satisfying this property. This finishes the proof. \qed

The next proposition is very similar to \cite{Ro} proposition 2.3.
However, there are some differences: the result of Rousseau holds in
a much more general situation; our result is formulated in terms of
the affine Weyl group, rather than in terms of the building, is valid only
for type $A_n$, and has an explicit proof, which gives some
quantitative information on the involved constants.

\begin{prop} \label{prop:LinFunctBeschr}
There is a linear function $f \colon \bZ_{> 0} \rar \bZ$, $f(z) = az + b$ depending only on $n$, with $a > 0$, such that for $v \in W_a$, we have $\ell(\beta(v)v^{-1}) \geq f(\ell(v))$. In particular, for a given $w \in W_a$, there are at most finitely many $v \in W_a$, such that $\beta(v)v^{-1} = w$. For any $r > 0$, the set
\[ \{ v \in W_a \colon \ell(\beta(v)v^{-1}) < r \} \]
\noindent is finite. Moreover, one can take $f$ to be
\[ f(z) = \frac{2}{n-1} z - \frac{2n(2n-3)}{n-1}. \]
\end{prop}

\begin{proof}
The second statement follows from the first: assume $w$ is given. Then $\beta(v)v^{-1} = w$ implies $\ell(w) > f(\ell(v))$. Since $a > 0$, there are at most finitely many positive integers $z$ with $f(z) < \ell(w)$. Thus the length of $v$ is bounded from above and hence there are at most finitely many such $v$'s. In particular, the finiteness of $\{ v \in W_a \colon \ell(\beta(v)v^{-1}) < r \}$ follows from the preceding statement.

Now we prove the first statement. Write $v = v_f \epsilon^{\lambda}$, with $v_f \in W$, $\lambda \in X_{\ast}(A)$. Then $\ell(v_f) \leq 2(n-1) - 1 = 2n - 3 =: c$ and therefore
\begin{equation}\label{eq:laengeofbetavv} \ell(\beta(v)v^{-1}) = \ell( \beta(v_f) \beta(\epsilon^{\lambda}) \epsilon^{-\lambda} v_f^{-1} ) \geq \ell( \beta(\epsilon^{\lambda}) \epsilon^{-\lambda} ) - 2c.
\end{equation}
\noindent If now $f(z) = az + b$ is a function, which has the property $\ell(\beta(\epsilon^{\lambda})\epsilon^{-\lambda}) \geq f(\ell(\epsilon^{\lambda}))$ for all $\lambda \in X_{\ast}(A)$, then $\tilde{f}(z) = az + b - (2 + a)c$ has this property for all $v \in W_a$. Indeed, for $v = v_f \epsilon^{\lambda}$, we have:
\begin{eqnarray*} \tilde{f}(\ell(v)) &=& a\ell(v) + b - (2 + a)c \leq a(\ell(\epsilon^{\lambda}) + c) + b - (2+a)c \\
&=& a\ell(\epsilon^{\lambda}) + b - 2c = f(\ell(\epsilon^{\lambda})) - 2c \leq \ell(\beta(\epsilon^{\lambda}) \epsilon^{-\lambda}) - 2c \leq \ell(\beta(v)v^{-1}),
\end{eqnarray*}
\noindent where the last inequality follows from \eqref{eq:laengeofbetavv}. Thus it is enough to prove the existence of the function $f$ satisfying the announced property only for elements $v = \epsilon^{\lambda}$ with $\lambda \in X_{\ast}(A)$. Write $\lambda = ( \lambda_1, \lambda_2, \dots, \lambda_n )$. Recall that for any $i \in \bZ$, we set $\lambda_i := \lambda_j$, where $j$ is unique with $1 \leq j \leq n$ and $i \equiv j \mod n$. For any $1 \leq k \leq n-1$, define
\[S_k := \sum\limits_{1 \leq i \leq n} \abs{\lambda_{i+k} - \lambda_i}. \]
\noindent Then $\sum_{k=1}^{n-1} S_k = \sum\limits_{1 \leq i \neq j \leq n} \abs{\lambda_j - \lambda_i}$ and an easy computation shows:
\[ 2\ell(v) = \sum_{1 \leq i \neq j  \leq n} \abs{\lambda_i-\lambda_j} = \sum_{k=1}^{n-1} S_k.\]
\noindent Further $\beta(v)v^{-1} = \epsilon^{\mu}$, where $\mu = (\lambda_{1+m}-\lambda_1, \lambda_{2 + m} - \lambda_2 , \dots, \lambda_m - \lambda_n)$. Thus we have
\[ 2\ell(\beta(v)v^{-1}) = \sum_{1 \leq i \neq j  \leq n} \abs{\lambda_i-\lambda_{i+m} - \lambda_j + \lambda_{j+m}}.\]
\noindent For $k \in \bZ$, let $0 \leq d(m,k) < n$ be such that $k \equiv md(m,k) \mod n$. We have:

\begin{lm} $S_k \leq d(m,k) S_m$
\end{lm}
\begin{proof}[Proof of the lemma] Write $d := d(m,k)$. Then
\[ dS_m = d \sum\limits_{i=1}^n \abs{\lambda_{i+m} - \lambda_i} = \sum\limits_{j = 1}^n \sum\limits_{l = 0}^{d-1} \abs{\lambda_{j + (l+1)m} - \lambda_{j + lm}}, \]
where the last equality is obtained by rearranging the terms and using that $(m,n) =  1$.
\noindent Then by triangle inequality we obtain: $\sum\limits_{l = 0}^{d-1} \abs{\lambda_{j + (l+1)m} - \lambda_{j + lm}} \geq \abs{\lambda_{j+dm} - \lambda_j} = \abs{\lambda_{j+k} - \lambda_j}$. This proves the lemma.
\end{proof}

Again by triangle inequality, $n \abs{\lambda_i - \lambda_{i+m}} \leq \sum_{j=1}^n \abs{(\lambda_i - \lambda_{i+m} - \lambda_j + \lambda_{j+m})}$. Summed over all $i$ between $1$ and $n$, this implies $n S_m \leq 2\ell(\beta(v)v^{-1})$. Putting together, we obtain:
\[ 2\ell(v) = \sum_{k=1}^{n-1} S_k \leq \sum_{k=1}^{n-1} d(m,k) S_m = \frac{(n-1)n}{2} S_m \leq  (n-1) \ell(\beta(v)v^{-1}),\]
since if $k$ runs through integers between 1 and $n-1$, then $d(m,k)$ also runs through integers between $1$ and $n-1$. This finishes the proof.
\end{proof}

\renewcommand{\refname}{References}

\end{document}